\documentstyle[12pt]{article}
\newlength{\abstractwidth}
\renewcommand{\thefootnote}{\fnsymbol{footnote}}
\renewcommand{\thanks}[1]{\footnote{#1}} 
\newcommand{\starttext}{
\setcounter{footnote}{0}
\renewcommand{\thefootnote}{\arabic{footnote}}}

\newcommand{\be}{\begin{equation}}
\newcommand{\bea}{\begin{eqnarray}}
\newcommand{\eea}{\end{eqnarray}}
\newcommand{\ee}{\end{equation}}

\def\ba{\begin{eqnarray}}
\def\ea{\end{eqnarray}}

\def\v{\vskip .1in}

\def\al{\alpha}

\def\d{\delta}
\def\e{\epsilon}
\def\g{\gamma}
\def\l{\lambda}
\def\m{\mu}
\def\n{\nu}
\def\o{\omega}

\def\r{\rho}

\def\t{\theta}

\def\D{\Delta}
\def\O{\Omega}

\def\cB{{\cal B}}

\def\cD{{\cal D}}

\def\cL{{\cal L}}

\def\cO{{\cal O}}

\def\cX{{\cal X}}

\def\Z{{\bf Z}}

\def\R{{\bf R}}
\def\C{{\bf C}}
\def\P{{\bf P}}
\def\K{{\rm K\"ahler }}

\def\Aut{{\rm Aut}}

\def\ti\tilde

\def\u{\underline}

\def\i{\infty}
\def\I{\int}
\def\p{\prod}
\def\s{\sum}

\def\ddb{\partial\bar\partial}
\def\sub{\subseteq}
\def\ra{\rightarrow}

\def\LA{\langle}
\def\RA{\rangle}

\def\det{{\rm det}}

\def\ti{\tilde}

\def\o{\omega}

\def\[{{\bf [}}
\def\]{{\bf ]}}

\v

\begin{document} \starttext \baselineskip=15pt
\setcounter{footnote}{0} \newtheorem{theorem}{Theorem}
\newtheorem{lemma}{Lemma}
\newtheorem{proposition}{Proposition}
\newtheorem{corollary}{Corollary}
\newtheorem{definition}{Definition}
 \begin{center} {\Large \bf DELIGNE
PAIRINGS AND THE KNUDSEN-MUMFORD EXPANSION}
\footnote{Research
supported in part by National Science Foundation grants DMS-02-45371
and DMS-05-14003}
\\
\bigskip

{\large D.H. Phong$^*$, Julius Ross$^*$ and
Jacob Sturm$^{\dagger}$} \\

\bigskip

$^*$ Department of Mathematics\\
Columbia University, New York, NY 10027\\

\v

$^{\dagger}$ Department of Mathematics \\
Rutgers University, Newark, NJ 07102

\end{center}
\v
\v\v\v \begin{abstract}

{\small Let $X\ra B$ be a proper flat morphism between
smooth quasi-projective varieties of
relative dimension $n$, and
$L\ra X$ a line bundle which is
ample on the fibers.
We establish formulas for the first two terms in
the Knudsen-Mumford expansion for
$\det (\pi_* L^k)$ in terms of  Deligne pairings
of $L$ and the relative canonical bundle~$K$.
This generalizes the theorem of Deligne \cite{Del} which holds
for families of relative dimension~one.
As a corollary,
we show that when $X$ is smooth, the line bundle
$\eta$
associated to $X\ra B$, which was
introduced in Phong-Sturm \cite{PS02}, coincides
with the CM bundle  defined by Paul-Tian \cite{PT06,PT06II}.
In a second and third corollaries, we
establish asymptotics for the K-energy
along Bergman rays, generalizing the
formulas obtained in \cite{PT06II}}.
\end{abstract}

\bigskip
\baselineskip=15pt
\setcounter{equation}{0}
\setcounter{footnote}{0}

\section{Introduction}
\setcounter{equation}{0}

Let $\pi: X\ra B$ be a flat proper  morphism of integral schemes
with constant relative dimension~$n$, and let
$ L\ra X$ be a relatively ample line bundle.
The theorem of Knudsen-Mumford \cite{KM} says that
there exist functorially defined line bundles
$\l_j=\l_j( X,L,B)\ra B$ with the
property:

\be\label{MK}
\det\ \pi_*(L^k)\ \approx \
\l_{n+1}^{\left(k\atop n+1\right)}
\otimes
\l_{n}^{\left(k\atop n\right)}
\otimes\cdots \otimes \l_0\ \ \
\hbox{
for\ $k>>0$
}
\ee
\def\cY{{\cal Y}}

In the case $n=1$, Deligne \cite{Del} showed that $\l_2(L,X,B)=\LA L,L\RA_{X/B}$,
the Deligne pairing of $L$ with itself.
If in addition the varieties $X$ and $ B$ are smooth, Deligne proved  that
$\l_1(L,X,B)^2=\LA LK^{-1},L\RA_{X/B}$, where $K=K_{X/B}=K_X\otimes K_B^{-1}$ is the
relative canonical line bundle.
Our first result provides a generalization of these formulas to the case where $n\geq 0$:

\begin{theorem}\label{main}
Let $\pi:X\ra B$ be a proper flat morphism of integral schemes
of  relative dimension~$n\geq 0$ and let
$L\ra X$ be a  line bundle which is very ample on the fibers.

\begin{enumerate}
\item There is a canonical functorial isomorphism
\be\label{iso1}
\l_{n+1}(L,X,B) \ = \ \langle L,...,L\rangle_{X/B}\
\ee

\item
If $X$ and $B$ are smooth, and $K$ is the relative canonical line bundle of $X\ra B$, then
there is a canonical functorial isomorphism
\be\label{iso3}
\l_n^2(L,X,B)\ = \ {\langle L^nK^{-1},L,...L\rangle_{X/B} }
\ee

\end{enumerate}
where the right sides of (\ref{iso1}) and (\ref{iso3}) are Deligne pairings of $n+1$ line bundles.
\end{theorem}

{\it Remark 1.} Knudsen-Mumford prove that the leading term, $\l_{n+1}(L,X,B)$,
is equal to the Chow bundle, so (\ref{iso1}) follows immediately by combining their
result with
the theorem of Zhang~\cite{Z96}, which gives a formula for
the Chow bundle in terms of the Deligne pairing
$\LA L,...,L\RA$. Thus the main content of
Theorem \ref{main} is the isomorphism (\ref{iso3}).

\v

{\it Remark 2.} We need to make precise the meaning of ``functorial''
in statements (\ref{MK}), (\ref{iso1}) and (\ref{iso3}).
This will be done in $\S 4$.

\v

{\it Remark 3.} The proof of Theorem \ref{main} generalizes to the
case where $X\ra B$ is a relative
complete intersection, but for simplicity of exposition, we
restrict to the smooth setting.
\v
{\it Remark 4.} If the base $B$ is compact, then it
is easy to see that both sides of (\ref{iso1}) and (\ref{iso3})
have the same Chern class. But in the statement of
Theorem \ref{main}, the base $B$
is not assumed to be compact. This means that we can not use a
Chern class argument and that we must work directly with
the sections of the relevant line bundles.  Allowing $B$ to
be non-compact is important for applications
to \K geometry, where cases of interest are test configurations $X\ra B$
for which the base $B$ is the complex plane.
\v
We shall give proofs of (\ref{MK}), (\ref{iso1}) and (\ref{iso3}) which
use Bertini's theorem and
go by induction on the dimension~$n$. The proofs are
self-contained, and do not rely on the $n=1$ result of
Deligne, or the results of Knudsen-Mumford \cite{KM}
and Zhang \cite{Z96}.
\v
Next we describe two applications of
Theorem \ref{main}. The first says that
the  line bundle
$\eta$, which was introduced
in \cite{PS02}, coincides with the CM line bundle
$\eta_{{\rm CM}}$, which was  defined by
Paul-Tian \cite{PT06, PT06II}.
In order to state the precise result, we first
recall the necessary definitions.

\v
Let $\pi:X\ra B$ be
a flat proper map of smooth quasi-projective varieties of
relative dimension $n$, and
let $L\ra X$ be a line bundle
which is relatively ample on the
fibers. Let $p(k)=a_0k^n+a_1k^{n-1}+\cdots $ be the
Hilbert polynomial of the fibers of $\pi$, and
let $\m={2a_1\over a_0}={nc_1(X)\cdot c_1(L)^{n-1}\over c_1(L)^n}$.
Let $K=K_{X/B}$ be
the relative canonical bundle of $X\ra B$.
\v
The line bundle $\eta\ra B$ was introduced
in \cite{PS02}
and it is defined as follows:

\be\label{PS}
\eta(L,X)\ = \ \langle L,...,L\rangle^\m
\otimes
\langle K,L,...,L\rangle^{(n+1)}\ ,
\ee
where, in each of the two Deligne pairings, there is
a total of $n+1$ line bundles. In \cite{PS02} it is  proved
that the
line bundle $\eta$ has a metric given by the Mabuchi K-energy.
In \cite{PS03} this
point of view was developed and
applied to the calculation of the
K-energy of a complete intersection, thus
providing a non-linear generalization of the
Futaki invariant formulas discovered by Lu
\cite{Lu99} and Yotov \cite{Y99} (see also the related
work of Hou \cite{Hou}, Liu \cite{Liu}, and Yotov \cite{Y98}).

\v
{\it Remark 5.} We learned recently that the line bundle $\eta$
had also been deduced from Riemann-Roch by Shou-Wu Zhang,
in a 1993 letter to P. Deligne \cite{Z93}.

\v

The
line bundle $\eta_{\rm CM}\ra B$ was introduced in
\cite{PT06} and it is defined
as follows:
\be \label{PT}
\eta_{{\rm CM}}(L,X)\ = \ \l_{n+1}(L,X,B)^\m\otimes
\left(
{\l_{n+1}(L,X,B)^n\over \l_n(L,X,B)^2}
\right)^{(n+1)}\ .
\ee

\v

This extends to a bundle on the Hilbert scheme whose weights are
the Donaldson-Futaki invariants which were defined in \cite{D02}.
\v
\begin{corollary}\label{cor one}
Let $\pi:X\ra B$ be
a flat proper map of smooth quasi-projective varieties and
$L\ra X$ a line bundle which is very ample on the fibers. Then
there is a canonical functorial isomorphism
$\eta(L,X) \ra \eta_{\rm CM}(L,X)$.

\end{corollary}

Before stating the next corollary we need
to recall some background from \K geometry
(full definitions will be provided in \S 6):
Let $X_1$ be a compact complex manifold and
$ L_1\ra X_1$ an ample line bundle. Then
Donaldson \cite{D02} defines a test
configuration $T$ for $(X_1,L_1)$ to be a triple
$L\ra X\ra \C$ plus a homomorphism
$\r: \C^\times~\ra \Aut(L,X,\C)$, where $X$ is a scheme,
 $L\ra X$ is a $\C^\times$
equivariant line bundle, $\pi:X\ra\C$ is flat and
$\C^\times$ equivariant,  $\pi^{-1}(1)=X_1$
and $L|_{X_1}= L_1$. Donaldson associates to $T$ a rational
number $F(T)$ which is called the Futaki invariant
of $T$,
and which
is defined by
$$
\r(\tau)|_{\eta_{\rm CM}(L,X)_0}=
\tau^{-F(T)}
$$
where
${\eta_{\rm CM}(L,X)_0}$ is the
fiber of ${\eta_{\rm CM}(L,X)}\ra\C$
at the origin.
This invariant generalizes the invariant
defined by Tian \cite{T97} in the case where the
central fiber is normal.
We say that $(X_1,L_1)$ is K-stable if $F(T)\le 0$ for all test
configurations $T$, with equality if and only if $T$ is a product.
The  conjecture of Yau \cite{Y87}, Tian \cite{T97}
and Donaldson \cite{D02}
 says that $X_1$ has a metric of constant
scalar curvature in $c_1(L_1)$ if and only if
the pair
$(X_1,L_1)$ is K-stable.

\v

Now assume that $L\ra X\ra\C$ is a test configuration
with $X$ smooth. Let $\o$ be a \K metric on $X$, and
let $\o_t=\o|_{_{X_t}}$ where for $t\in\C^\times$, we
put $X_t=\pi^{-1}(t)$. Let $d=\I_{X_1}\o_1^n$ and define, as in \cite{T97}, the function

\be \psi(t)\ = \ {1\over d}\I_{X_t} \log
\left[{\o^n\wedge (d\pi\wedge d\bar \pi)\over \o^{n+1}}\right]\ \o_t^n
\ee
Here we view $\pi:X\ra {\bf C}$ as a holomorphic function,
so that $d\pi$ is a 1-form on $X$. The expression
$f={\o^n\wedge (d\pi\wedge d\bar\pi)\over\o^{n+1}}$
is the ratio of two $(n+1,n+1)$ forms on $X$ whose denominator is strictly
positive, and thus $f$ is a non-negative smooth function on $X$.
Then $\psi:\C^\times\ra\R$ is smooth and bounded above as $t\ra 0$.
Let $\n(t)=\n(\o_1,\r(t)^*\o_t)$ be the K-energy of $\r(t)^*\o_t$
with respect to the base point $\o_1$.
Then we have the following generalization of the formula
proved in \cite{PT06II}:

\v

\begin{corollary}\label{asym}
 Let $L\ra X\ra\C$ be a test configuration with $X$ smooth,
and let $\o\in c_1(L)$ be a \K metric on $X$. Then
\be \n(t)\ -  \ \psi(t)\ = \ {F(T)\over (n+1)d}\log |t|^2\ + \ O(1)
\ee
Hence, if the central fiber of $X$ has no component of multiplicity
greater than one, then
\be \lim_{t\ra 0} \ {\n(t)\over \log |t|^2}\ = \ {F(T)\over (n+1)d}
\ee
\end{corollary}

This result was obtained in \cite{PT06II} under the following additional assumption:
There is a triple
$(\cL,\cX,\cB)$ with $\cX\ra \cB$ a flat map between smooth projective varieties, $\cL\ra\cX$ relatively
very ample, $\P(\pi_*\cL)\approx B\times \P^N$, $(X,L)\approx (\cX_b,\cO_{\P^N}(1))$ for some
$b\in\cB$, with the property: There is
an action of $SL(N+1,\C)$ on the data commuting with all the projections such
that $\r$ is the restriction of a one parameter subgroup of $SL(N+1,\C)$.
The purpose of Corollary \ref{asym} is to remove this assumption.

\v

More generally, suppose $L\ra X\ra \C$ is a test configuration $T$ and
$L'\ra X'\ra B$ a flat family satisfying the hypothesis of
 Corollary \ref{cor one}. Suppose that there is an imbedding $\C\sub B$
such that $T$ is the restriction of the flat family to $\C$. Thus $X'$ is
smooth, but $X$ need not be smooth. Then we have the following generalization
of Theorem 1 in \cite{PT06II} :

\begin{corollary}
Let $\o\in c_1(L')$ be a \K metric on $X'$ and $\eta$ a \K metric on $B$. Define
\be
 \psi(t)\ = \ {1\over d}\I_{X_t} \log
\left[{\o^n\wedge \pi^*\eta^m\over \o^{n+m}}\right]\ \o_t^n
\ee
where $n+m$ is the dimension of $X$. Then
\be \n(t)\ -  \ \psi(t)\ = \ {F(T)\over (n+1)d}\log |t|^2\ + \ O(1)
\ee
\end{corollary}
\v\v

{\bf Acknowledgement}. We would like to thank Shou-Wu Zhang for some useful conversations,
for a careful reading of the paper and many helpful suggestions,
and for informing us about his  1993 correspondence
with Pierre Deligne.

\section{Review of the Deligne pairing}

\setcounter{equation}{0}
We recall some of the results of Deligne \cite{Del} and Zhang \cite{Z96}:
Let $\pi:X\ra B$ be a flat projective morphism of integral schemes
of pure relative dimension~$n$. If $L_0,...,L_n$ are line bundles
on $X$, then the Deligne pairing $\langle L_0,...,L_n\rangle(X/B)$
is a line bundle on $B$. It is locally generated by symbols
$\langle s_0,...,s_n\rangle$, where the $s_j$ are rational sections
of the $L_j$ whose divisors, $(s_j)$, have empty intersection.
The transition functions are determined by the following relation:

\be\label{del def} \langle s_0,..., fs_i,...s_n\rangle \ = \ f[Y]\cdot
\langle s_0,...,s_n\rangle
\ee
where $f$ is a rational function, $Y= \cap_{j\not=i} (s_j)$ is flat
over $B$, and $f[Y]={\rm Norm}_{Y/B}(f)$. The fact that
(\ref{del def}) determines a well defined line bundle follows
from the Weil reciprocity formula, which says that if $f$ and $g$
are rational functions on a projective curve such that
$(f)\cap (g)=\emptyset$, then

\be\label{weil}
f[(g)]\ = \ g[(f)]
\ee
where $g[(f)]=\prod_{j} g(p_j)^{\m_j}$ if $(f)=\s_j \m_jp_j$.
\v
The Deligne pairing satisfies a useful induction formula which
is described as follows: Suppose $s_j$ is a rational section of
$L_j$ such that $Y\ra B$ is flat where $Y=(s_j)$. Then the map
$\langle s_0|_Y,..., \hat s_j|_Y,...,s_{n}|_Y\rangle(Y/B)  \ra \langle s_0,...,s_{n}\rangle(X/B)$
defines an isomorphism

\be\label{induct}
 \langle L_0|_Y,...,\hat L_j|_Y,...,L_{n}|_Y\rangle (Y/B) \ra \
\langle L_0,...,L_n\rangle(X/B)
\ee

Next let us suppose that for some $i<j$ that $L_i=L_j$,
and let $\langle s_0,...,s_n\rangle$ and $\langle t_0,...,t_n\rangle$
be generating sections for $\langle L_0,...,L_n\rangle$
as above. Assume that $s_k=t_k$ for all $k\not=i,j$. Assume
also that $s_i=t_j$ and $s_j=t_i$. Then we have the following
formula from \cite{Del}

\be\label{permute} \langle s_0,...,s_n\rangle
\ = \
(-1)^d\langle t_0,...,t_n\rangle
\ee

where $d=\p_{k\not=i}c_1(L_k)$.
In \cite{Del}, the relation (\ref{permute})
is only stated  for $n=1$, but the general
case follows from this and from (\ref{induct}).

\v

Finally, we recall that if $h_j$ is a smooth metric
on $L_j$, then there is an induced metric  $\langle h_0,...,h_n\rangle$
on the line bundle $\langle L_0,...,L_n\rangle$. This metric has
the following property: Let $\phi$ be a smooth function on $X$.
Then $h_0e^{-\phi}$ is a metric on $L_0$ and

\be \langle h_0e^{-\phi},...,h_n\rangle=\langle h_0,...,h_n\rangle\cdot e^{-\psi}
\ee
where $\psi: B\ra \C$ is the function

\be
\psi\ = \ \I_{X/B} \phi \cdot \o_1 \wedge\cdots \wedge\o_n
\ee

and $\o_j=-{\sqrt{-1}\over 2\pi}\ddb\log h_j$ is the curvature of $h_j$.

\v
\section{Bertini's Theorem}

Our proofs require a variant of Bertini's theorem
to cut down the dimension of our family $X\ra B$
so that we obtain a smooth family of smaller
relative dimension. We state the version of
the theorem that we need and, for the sake of
completeness, we supply a proof.

\v
\begin{proposition}\label{Bertini} Let $X\sub\P^N$
be a smooth quasi-projective subvariety and let
$Y\sub X$ be any subvariety with $Y\not=X$. Let $H_0\sub\P^N$ be a hyperplane
such that $H_0\cap X$ is smooth.  Let $\ell$ be a
generic pencil of hyperplanes containing $H_0$, that is, $\ell$ is
a generic line in $\P_*^N=\{H\sub\P^N: H$ is a hyperplane$\}$ which
passes through $H_0\in \P_*^N$. Then there exists an open
set $U\sub X$ such that $Y\sub U$ and $H\cap U$ is smooth
for all $H\in\ell$.
\end{proposition}

{\it Proof.} We first recall the proof of the usual Bertini theorem:
Let $X\sub\P^N$ be a smooth variety of dimension $n$ and let $\ell$
be a generic pencil of hyperplanes. We wish to show that $H\cap X$ is smooth
for all but finitely many $H\in\ell$. To see this, let
$F: X\ra X\times {\bf Gr}(n,N)$ (the set of $n$ planes in $\P^N$) be the
map which sends $x$ to $(x,T_x(X))$. Let $p: X\times {\bf Gr}(n,N)\ra {\bf Gr}(n,N)$ be
the projection map. Let $Z\sub {\bf Gr}(n,N)\times
\P^N_*$ be the set of pairs $(\l,H)$ such that $\l\sub H$. Then $\pi_1: Z\ra {\bf Gr}(n,N)$
has fibers of dimension $N-n-1$. Let $\pi_2:Z\ra {\bf P}_*^N$
be the projection onto the second factor.
Thus $\cB=(\pi_2\pi_1^{-1}pF)(X)\sub \P^N_*$ is a
constructible set of dimension at most $N-1$. Moreover, $H\in  \cB\iff H\cap X$
is not smooth. Thus most hyperplanes (a non-empty Zariski open set) will
intersect $X$ along a smooth divisor. On the other hand, a line in $\P^N_*$
(i.e., a pencil of hyperplanes)
will, in general, contain a finite number of  hyperplanes $H$ for
which $H\cap X$ is singular.

\v

Now we prove the proposition: Let $\cB_Y=(\pi_2\pi_1^{-1}pF)(Y)\sub \P^N_*$.
Then $\cB_Y$ is constructible and has dimension at most $N-2$. Thus
a generic pencil $\ell$ misses $\cB_Y$. Let $\{H_1,...,H_k\}=\{H\in \ell: H\cap X$ is not smooth $\}$.
Let $U_j\sub H_j\cap X$ be the set of smooth points. Then $Y\sub U_j$.
Then $U=\cap U_j$ is the desired open set, proving the proposition.
\v

Next we let
 $f:X\ra B$ be a projective flat morphism between smooth varieties and
$L\ra X$ a line bundle which is very ample on fibers. Assume that $B$ is
affine and fix $b_0\in B$.
Assume as well that $X\sub B\times \P^N$
is the imbedding given by the complete linear series of $L$.
\v

\begin{proposition}\label{Bertini 2}
There exists  $s\in H^0(X,L)$ such that $\{s=0\} \sub X$ is
smooth and flat over $B'$ where $b_0\in B'\sub B$ is
open. Moreover, if $s_1,s_2\in H^0(X,L)$ are two such
sections, then there exists $s'\in H^0(X,L)$ such that
$\{ts_1+(1-t)s'=0\}$ and $\{ts_2+(1-t)s'=0\}$ are smooth
and flat over $B''$ for all $t\in\C$ (where $b_0\sub B''\sub B'$
is open).
\end{proposition}

{\it Proof.}
Let $f:B\ra \C^m$ be an imbedding where
$f=(f_0,...,f_m)$ and $f_j:B\ra \C$ a regular function, and $f_0=1$.
Let $x_0,...,x_N$ be the homogeneous coordinates on $\P^N$.
Then the map $\P^N\times B\ra \P^M$ given by
$(x,b)\ra (f_ix_j)$ is an imbedding which restricts
to an imbedding of $X\hookrightarrow \P^M$. Here $M=(m+1)(N+1)-1$.
Thus, if
$c_{ij}$ are generic constants, the usual Bertini
theorem says that $s=\s_{ij}c_{ij}f_ix_j=0$
is a smooth subvariety of~$X$. And this subvariety
is flat over $B$ (possibly after shrinking $B$ a little).
\v
Now let $s_1,s_2\in H^0(X,L)$. Then there exist
regular functions $f_0,...,f_m$ with $f_0=1$ such that $s_i$ is
a $\C$ linear combination of the sections $f_ix_j$ and
such that $(x,b)\ra (f_ix_j)$ is an imbedding of $\P^N\times B$.
The existence of $s'$ now follows from Proposition \ref{Bertini}.

\v
\section{Formula for the leading term }
\setcounter{equation}{0}
In this section we prove (\ref{MK}) and (\ref{iso1}).
Although the results are not new, a methodology will be developed
which will also yield, with suitable adaptations,
a proof of  (\ref{iso3}).

\v

Before constructing the isomorphisms of (\ref{MK}) and (\ref{iso1}), we
must first make precise the meaning of
``functorial'':
Let $B$ be a scheme and suppose
we are given, for all sufficiently large positive integers $k$,  a line
bundle $M_k\ra B$. Then associated
to such a sequence $M^\bullet=(M_k)$ we
define, as in \cite{KM}, a new sequence $\D M^\bullet $ as follows:
$(\D M^\bullet)_k= M_k\otimes M_{k-1}^{-1}$.

\v

Now
let $\pi: X\ra B$ be a flat
proper morphism of integral
schemes of relative dimension $n$, and $L\ra X$ a line bundle
which is relatively ample on the
fibers.
Then a precise formulation of the
Knudsen-Mumford theorem (\ref{MK}) is the following:

\begin{theorem}\label{reformulation}
There exists, for
$k>>0$, an isomorphism
\be\label{KN1}
\sigma_k(L,X,B): \D^{(n+2)}\det(\pi_*L^k)
\ \ra \ \cO_B
\ee
with the functorial property:
If $\phi:(L',X',B')\ra (L,X, B)$ is a cartesian
morphism
then
\be\label{functorial}
\sigma_k(L',X',B')\circ\phi^*=\phi^*\circ
\sigma_k(L,X,B)
\ee
where
$\phi^*$
denotes the maps $\cO_B\ra\cO_{B'}$
and $\D^{(n+2)}\det(\pi_*L^{\bullet})_k
\ra \D^{(n+2)}\det(\pi'_*L^{'\bullet})_k
 $
induced by $\phi$.
\end{theorem}

Here we say that $\phi=(\phi_2,\phi_1,\phi_0)$ is cartesian
if $\phi_0:B'\ra B$ and $\phi_1:X'\ra X$ are morphisms such
that $\pi\phi_1=\phi_0\pi'$, the induced map $X'\ra X\times_B B'$
is an isomorphism, and $\phi_2: L'\ra \phi_1^*L$ is an isomorphism
of line bundles over $B'$.

\v
Let us spell out the equivalence
between (\ref{KN1}) and (\ref{MK}), which is the usual
formulation  of the Knudsen-Mumford
theorem: Fix $n\geq 0$ and assume that
(\ref{KN1}) holds. For $k>>0$,
we define
\be\label{lamdanplusone}
\l_{n+1}(X,L,k)=\D^{(n+1)}
\det(\pi_*L^{\bullet})_{k}\ee
Then
(\ref{KN1}) defines a family of functorial
isomorphisms $\sigma_{k,k'}(n+1):
\l_{n+1}(X,L,k)
\ra
\l_{n+1}(X,L,k')$. Now fix $k_0>>0$
and define $\l_{n+1}(X,L)=\l_{n+1}(X,L,k_0)$.
Let
\be\l_{n}(X,L,k)=\ \D^{(n)}
\left\{
\det(\pi_*L^\bullet)\otimes \left[\l_{n+1}(X,L)^{\left({\bullet\atop n+1}\right)}\right]^{-1}\right\}_k
\ee
Then $\sigma_{k,k'}(n+1)$ defines
functorial
isomorphisms $\sigma_{k,k'}(n):
\l_{n}(X,L,k)
\ra
\l_{n}(X,L,k')$.
Let $\l_n(X,L)=\l_n(X,L,k_0)$ and define

$$ \l_{n-1}(X,L,k)\ = \
\D^{(n-1)}
\left\{
\det(\pi_*L^\bullet)\otimes \left[\l_{n+1}(X,L)^{\left({\bullet\atop n+1}\right)}\right]^{-1}
\otimes
\left[\l_{n}(X,L)^{\left({\bullet\atop n}\right)}\right]^{-1}\right\}_k
$$

Continuing in this fashion, we
construct line bundles $\l_j(X,L)$
for $0\leq j\leq n+1$ which satisfy
(\ref{MK}). Moreover, if $\phi:B'\ra B$
is a base change, then the construction
of the $\l_j$ provides an isomorphism
$\phi^*\l_j(X,L)\ra \l_j(X',L')$ which
is compatible with the isomorphism
$\phi^*\det(\pi_*L^k)\ra
\det(\pi'_*{L'}^k)$.

\v
Next we give a precise formulation of the
statement  (\ref{iso1}) which gives the
formula for $\l_{n+1}(X,L)$
in terms of
a Deligne pairing:
\begin{theorem}\label{sharp}
Let $\pi: X\ra B$ be a flat
projective morphism of quasi-projective
schemes of relative dimension $n$, and $L\ra X$ a line bundle
which is very ample on the
fibers.
There exists, for
$k>>0$, an isomorphism
\be\label{KN}
\tau_k(L,X,B): \D^{(n+1)}\det(\pi_*L^{\bullet})_k
\ \ra \ \langle L,...,L\rangle
\ee
with the functorial property:
If $\phi:(L',X',B')\ra (L,X, B)$ is a cartesian
morphism then
\be\label{functorial}
\tau_k(L',X',B')\circ\phi^*=\phi^*\circ
\tau_k(L,X,B)
\ee
where
$\phi^*$
denotes the isomorphism $\phi^*\langle L,...,L\rangle
\ra\langle L',...,L'\rangle
$
as well as the isomorphism $\phi^*\D^{(n+2)}\det(\pi_*L^k)
\ra \D^{(n+2)}\det(\pi'_*{L'}^k)
 $
induced by $\phi$. In particular,
there is a functorial isomorphism
\be \l_{n+1}(L,X,B)\ \ra \
\langle L,...,L\rangle
\ee
\end{theorem}

Note that Theorem \ref{sharp} implies
Theorem \ref{reformulation} upon
defining $\sigma_\bullet(L,X,B)=\D\tau_\bullet(L,X,B)$.
\v
Finally we give the precise formulation of (\ref{iso3}):

\begin{theorem}\label{sharp1}
Let $\pi: X\ra B$ be a flat
projective morphism of smooth quasi-projective
varieties of relative dimension $n$, and $L\ra X$ a line bundle
which is very ample on the
fibers.
There exists, for
$k>>0$, an isomorphism
\be\label{phik}
\m_k=\m_k(L,X,B): [\D^{(n)}\det \pi_*(L^k)]^2
\ \ra \
\langle L,...L\rangle^{2k}\langle K_{X/B}L^n,...,L\rangle^{-1}
\ee
with the functorial property:
If $\phi:(L',X',B')\ra (L,X, B)$ is a cartesian
morphism then
\be\label{functorial}
\m_k(L',X',B')\circ\phi^*=\phi^*\circ
\m_k(L,X,B)
\ee
In particular,
there is a functorial isomorphism
\be
\l_n^2(L,X,B)\ = \ {\langle L^nK^{-1},L,...L\rangle_{X/B} }
\ee
\end{theorem}
\v

{\it Proof of Theorem \ref{sharp}}: We will
define $\tau_k(L,B)$ first on the level of
stalks. Thus, we
fix $b_0\in B$, and we define $\tau(L,B')$ where $b_0\in B'\sub B$
is some small open neighborhood. In order
to avoid cumbersome notation, we shall write $B$
instead of $B'$, with the understanding that
$B$ has possibly been replaced by a smaller open neighborhood of $b_0$.
The definition we give
will depend on some choices, so the main task
will be to show that after shrinking $B$ even
further, that different sets of choices
define the same $\tau_k(L,X,B)$.

\v
We start with the case of relative dimension zero:
Fix a section $s$ which generates $L$ (shrinking the base $B$
if necessary).
Then multiplication by $s$ defines an isomorphism
between $L^{k-1}$ and $L^k$ and thus an isomorphism
$\det(\pi_*L^{k-1})\ra \det( \pi_*L^{k})$. This
provides a nowhere vanishing section $\xi$ of $\D\det(\pi_*L^k) $
and the map $\langle s\rangle \mapsto \xi$ defines
$\tau_k(L,B)^{-1}$.

\v
Now let $n$, the relative dimension, be arbitrary.
 Choose generic sections
$s_1,...,s_n$ of $\pi_*L$ and for $0\leq k\leq n$,
let $X_k\sub X$ be the subscheme defined
by $s_{k+1}=s_{k+2}=\cdots =s_n=0$. Thus $X_n=X$
and, applying Proposition \ref{Bertini 2}, we conclude that $X_j\ra B$
is a projective flat map between smooth
quasi-projective varieties (again, with  the
understanding that
$B$ has been replaced by a smaller open neighborhood of~$b_0$).
\v

Multiplication by $s_j$ defines an exact sequence

\be
0\ \ra \ \pi_*L^{k-1}_{X_j}\ \ra \ \pi_*L^{k}_{X_j}\ \ra \
\pi_*L^{k}_{X_{j-1}}\ \ra \ 0
\ee

Taking determinants defines an isomorphism

\be
\kappa_{s_j}: \D \det(\pi_*L_{X_j}^k) \ \ra \ \det(\pi_*L|_{X_{j-1}}^k)
\ee
More precisely, if $t_1,...,t_a$ is a basis of $H^0(B,\pi_*L^{k-1}_{X_j})$, choose
$u_1,...,u_b\in H^0(B,\pi_*L^k_{X_j})$ so that
$\{s_jt_1,...,s_jt_a, u_1,...,u_b\}$ is a basis of $H^0(B,\pi_*L^k_{X_j})$. Then
we define

$$\kappa_{s_j}\left(( s_jt_1\wedge\cdots\wedge s_jt_a\wedge u_1\cdots \wedge u_b)\otimes
(t_1\wedge \cdots \wedge t_a)^{-1}\right)\ = \ \tilde u_1\wedge\cdots \wedge\tilde u_b
$$
where $\tilde u_i=u_i|X_{j-1}$.

\v

Define the isomorphism
$\kappa(s_1,...,s_n)=\D\kappa_{s_1}\circ \D^{(2)}\kappa_{s_2}\circ\cdots \circ
\D^{(n)}\kappa_{s_n}$. Then

\be \kappa(s_1,...,s_n): \D^{(n+1)}\det(\pi_*L^k)\ \ra \ \D\det(\pi_*L^k_{X_0})
\ee

\v
Now let $\iota_{s_j}:\langle L|_{X_j},...,L|_{X_j}\rangle_{X_j/B}
\ra \langle L|_{X_{j-1}},...,L|_{X_{j-1}}\rangle_{X_{j-1}/B}$
be the induction isomorphism defined by (\ref{induct}).
Define the isomorphism $\iota(s_1,...,s_n)=\iota_{s_1}\circ \cdots \circ \iota_{s_n}$. Then

\be \iota(s_1,...,s_n): \langle L,...,L\rangle_{X/B}\ \ra
\ \langle L|_{X_0}\rangle_{X_0/B}
\ee

Finally, define

\be
\tau_k(L,X,B)(s_1,...,s_n)\ = \ \iota(s_1,...,s_n)^{-1}\circ \tau_k(L|_{X_0},X_0,B)\circ\kappa(s_1,...,s_n)
\ee

\begin{lemma}\label{ind}
The isomorphism $\tau_k(L,X,B)(s_1,...,s_n)$ is independent of the choice
of generic sections $s_1,...,s_n$.
\end{lemma}

Once Lemma \ref{ind} is proved, we can define $\tau_k(L,X,B)=\tau_k(L,X,B)(s_1,...,s_n)$
for any choice of defining sections. This is a local definition, since $B$ has
been replaced by a small open subset of $b_0$. Now Lemma \ref{ind} implies
that the local isomorphisms glue together to give a unique isomorphism which
is globally defined and easily seen to satisfy the functorial properties:
If  $s_1,\ldots,s_n$ are general elements of $H^0(L)$ used to define $\tau_k(X,L,B)$ then
we can use the pullbacks of these to $H^0(L')$, which are also general, to define
$\tau(X',L',B')$.
 Thus,
to prove Theorem \ref{sharp}, it suffices to prove the lemma.
\v
{\it Proof of Lemma 1.} First we show that $\tau_k$ is independent of the
ordering of the sections $s_1,...,s_n$. To see this, it suffices to show
that $\tau_k$ is invariant under a permutation which switches $s_{j-1}$ and $s_{j}$
for some $j<n$,
and fixes all the other sections. Let's verify this for $j=n$ (the general case
is similar): Let $\ti\kappa(s_{n-1},s_n)= \kappa_{s_{n-1}}\circ \D\kappa_{s_n} $. Thus
we have
$\ti\kappa(s_{n-1},s_n): \D^{(2)}\det(\pi_*L^k)\ra \det(\pi_*L^k|_{X_{n-2}})$.
From the definition of $\kappa_{s_j}$, we easily see that

\be \ti\kappa(s_n,s_{n-1})=(-1)^{p_{n-1}(k-1)}\ti\kappa(s_{n-1},s_n)
\ee
where $p=p_n$ is the Hilbert polynomial for $L\ra X\ra B$ and
$p_{j-1}(k)=\D p_{j}(k)=p_j(k)-p_j(k-1)$.
Indeed,
let $E_{k-1}$ be an ordered basis of $H^0(X,L^{k-1})$ . Let
$F_k\sub H^0(X,L^k)$ be an ordered set of sections whose restrictions to $X_{n-1}$
form a basis of $H^0(X_{n-1},L^k)$. Then $(s_{n-1}E_{k-1},F_k)$ is an
ordered basis of $H^0(X,L^k)$ and

$$\kappa_{s_n}\left({\det(s_nE_{k-1},F_k)\over \det(E_{k-1})}\right)\ = \ \det(F_k|_{X_{n-1}})
$$

We repeat this process with $X$ replaced by $X_{n-1}$: Let $F_{k-1}\sub H^0(X,L^{k-1})$ be
an ordered set such that $F_{k-1}|_{X_{n-1}}\sub H^0(X_{n-1},L^{k-1})$ is a basis.
Let $H_k\sub H^0(X,L^k)$ an ordered set such that
$(s_{n-1}F_{k-1}|_{X_{n-1}},H_k|_{X_{n-1}})$ is a basis of $H^0(X_{n-1},L^{k})$. So, starting with $F_{k-1}$ and
$H_k$ we define $F_k=(s_{n-1}F_{k-1},H_k)$. If, in addition, we fix a basis $E_{k-2}\sub H^0(X,L^{k-2})$,
we can define $E_{k-1}=(s_{n-1}E_{k-2},F_{k-1})$.
Then  we obtain

$$\kappa_{s_n}\left(
{\det(s_ns_{n-1}E_{k-2},s_nF_{k-1}, s_{n-1}F_{k-1},H_k)\over
\det(s_{n-1}E_{k-2},F_{k-1})}\right) \ =  \
\det(s_{n-1}F_{k-1}|_{X_{n-1}}, H_k|_{X_{n-1}})
$$

$$\kappa_{s_n}\left({\det(s_nE_{k-2},F_{k-1})\over \det(E_{k-2})}\right)\ = \ \det(F_{k-1}|_{X_{n-1}})
$$

Thus

$$\D\kappa_{s_n}
\left(
{\det(s_ns_{n-1}E_{k-2},s_nF_{k-1}, s_{n-1}F_{k-1},H_k)\det(E_{k-2})\over
\det(s_{n-1}E_{k-2},F_{k-1})\det(s_nE_{k-2},F_{k-1})}\right) \ =  \
{\det(s_{n-1}F_{k-1}, H_k)\over \det(F_{k-1})}
$$

where, to ease the notation, we've omitted the restrictions to $X_{n-1}$.
Applying $\kappa_{s_{n-1}}$ to both sides we obtain

$$ \tilde\kappa(s_{n-1},s_n)
\left(
{\det(s_ns_{n-1}E_{k-2},s_nF_{k-1}, s_{n-1}F_{k-1},H_k)\det(E_{k-2})\over
\det(s_{n-1}E_{k-2},F_{k-1})\det(s_nE_{k-2},F_{k-1})}\right) \ =  \ \det(H_k)
$$

If we interchange $s_n$ and $s_{n-1}$ in the ratio

$$\left(
{\det(s_ns_{n-1}E_{k-2},s_nF_{k-1}, s_{n-1}F_{k-1},H_k)\det(E_{k-2})\over
\det(s_{n-1}E_{k-2},F_{k-1})\det(s_nE_{k-2},F_{k-1})}\right)
$$

the sign changes by a factor of $(-1)^{|F_{k-1}|}$. On the other hand, (4.11)
implies that $|F_{k-1}|=p_{n-1}(k-1)$. This proves (4.16).

\v
Next, applying $\D$ successively, we have
$\D^{(n-1)}\ti\kappa(s_{n},s_{n-1})=(-1)^d\D^{(n-1)}\ti\kappa(s_{n-1},s_{n})$
where $d=\D^{(n-1)}p_{n-1}(k-n)=p_0(k-n)$. But $p_0=c_1(L)^n$ is a constant.
Thus $d=c_1(L)^n$. On the other hand, (\ref{permute}) implies that
permuting $s_n$ and $s_{n-1}$ in $\iota(s_1,...,s_n)$ introduces the same
factor of $(-1)^d$. Thus the two factors cancel, and we see that
$\tau_k$ is independent of the ordering of the $s_j$.
\v
Now let us fix two choices of sections, $(s_1,...,s_n)$ and $(s_1',...,s_n')$.
We must show that $\tau_k(L,X,B)(s_1,...,s_n)=\tau_k(L,X,B)(s_1',...,s_n')$.
Clearly we may assume that $s_i=s_i'$ for all but one index $i$. Since
$\tau_k$ is independent of the ordering of the sections, we may assume
that $s_j=s_j'$ for all $j\geq 2$. In other words, in the proof of Lemma 1,
we may assume that $n=1$.
\v
If $n=1$, then

\be\label{two sections}
\tau_k^{-1}(s_1)(\langle s_1,s_0\rangle)\ = \ {\det(E)\det(s_0s_1E,s_1F,s_0F)\over
\det(s_0E,F)\det(s_1E,F)}
\ee
Here $E$ is any basis of $H^0(\pi_*L^{k-2})$ and $F\sub H^0(\pi_*L^{k-1})$ is
any set of linearly independent elements such that $(s_0E,F)$, $(s_1E,F)$ and $(s_1'E,F)$
each form a basis of $H^0(\pi_*L^{k-1})$. We must show that
$\tau_k(s_1)(\langle s_1,s_0\rangle)\ = \tau_k(s_1')(\langle s_1,s_0\rangle)$.
Since $\langle s_1',s_0\rangle= f((s_0))\langle s_1,s_0\rangle$ where,
$f=s_1'/s_1$ and $(s_0)$ is the divisor of $s_0$, it suffices to show that
$\tau_k(s_1')(\langle s_1',s_0\rangle)=f((s_0))\tau_k(s_1)(\langle s_1,s_0\rangle)$.
But this follows immediately from (\ref{two sections}). This completes the
proof of Lemma \ref{ind} as well as the proofs of Theorem \ref{sharp}
and Theorem \ref{reformulation}.

\section{Relative dimension zero}

\setcounter{equation}{0}

In this section we prove Theorem \ref{sharp1} in the case $n=0$.

\v

\begin{lemma}\label{level one}
Let $\pi:Y\ra B$ be a finite flat morphism of integral schemes
and let $L\ra Y$ be a line bundle. Then
there is a canonical isomorphism

\be \det (\pi_*L)^{ 2} \ = \ \langle L\rangle^2\otimes \d_{Y/B}
\ee
where $\d_{Y/B}\sub\cO_B$ is the discriminant of the extension $Y\ra B$.
\end{lemma}

{\it Proof.} This question is local on the base, so we may assume $Y\ra B$ is a
finite morphism of degree $r$ between affine varieties. Thus $Y={\rm spec}(T)$ and $B={\rm spec}(S)$
where $S\sub T$ is a  extension of  rings such that $T$ is a free $S$
module of rank $r$.
\v

1) We must prove
\be\label{first} \det (\pi_*{\cal O})^2 = \d_{Y/B}
\ee

2) We must also show

\be\label{second} {(\det\ \pi_*L)\otimes (\det \ \pi_*\cO)^{-1}}\ = \ \langle L \rangle
\ee

To see (\ref{second}), let $s$ be a section of $L$. Then multiplication by $s$ defines
a map $\cO_Y\ra L$ and thus a map $\pi_*\cO\ra \pi_* L$ and hence a section of
${(\det\ \pi_*L)\otimes (\det \ \pi_*\cO_Y)^{-1}}$. One easily checks that
this defines the isomorphism asserted by (\ref{second}).
\v
As for (\ref{first}), let $\t_1,...,\t_r\in T$, where $r$ is the degree of $S\sub T$.
Let $\sigma_1,...,\sigma_r$ be the imbeddings of $T$ into the algebraic closure of
the fraction field of $S$. Then $\det(\sigma_i(\t_j))^2\in~\d_{Y/B}$. In fact, $\d_{Y/B}$ is
generated by all such elements, so this map gives the isomorphism (\ref{first})
and Lemma \ref{level one} is proved.

\v

\begin{lemma}\label{different}
Let $\pi:Y\ra B$ be a finite flat separable morphism of
quasiprojective varieties of relative dimension zero, and let
$\cD_{Y/B}\sub\cO_Y$ be the sheaf of ideals which annihilates
$\O_{Y/\cB}$, the module of relative differentials. Then $\cD_{Y/B}$
is locally free of rank one and
\be\label{6.1}\langle\cD_{Y/B}\rangle_{Y/B}=\d_{Y/B} \ee
\end{lemma}

{\it Proof.}
We may assume that $ B$ and $Y$ are affine:
$Y={\rm spec}(T)$ and $\cB={\rm spec}(S)$.
Then $\cD\sub T$ is an ideal.
Choose $\al\in T$ such that $T=S[\al]$. The
existence of such an $\al$ is guaranteed
by Nakayama's lemma (perhaps after shrinking the base).
Let $f$ be the minimal monic polynomial for $\al$,
so that $f(\al)=0$ and $\deg f=d$. Then, as is
well known,
\be\label{short}
\O_{Y/ B}\approx T/(f'(\al))
\ee

Let us briefly recall the proof. Since
$S[\al]=S[X]/(f)$ we have an exact sequence \cite{H}:

$$ (f)/(f)^2\ \ra \ \O_{S[X]/S}\otimes T\ \ra \ \O_{T/S}\ \ra 0
$$
where the first arrow is the map $u\mapsto du\otimes 1$.
The module in the middle is the free $T$
module generated by $dX\otimes 1$. The image of the first map is the free
$T$ module generated by $df \otimes 1 = f'(X)dX\otimes 1 = dX\otimes f'(\al)$.
This proves (\ref{short}).

\v

Now we can define the isomorphism of (\ref{6.1}): If $f'(\al)$ is a
generator of $\cD$ then we associate to it the basis
$\{1,\al,...,\al^{d-1}\}$ of $p_*\cO$ and thus an element
of $\det(p_*\cO)^2$. To see that this is well defined, let
$f'(\g)$ be another generator. Then, from the definition
of the Deligne pairing,
$$ \langle f'(\al)\rangle \ = \ {\rm Norm}_{T/S}
\left({f'(\al)\over f'(\gamma)}\right)\langle f'(\g)\rangle
$$

On the other hand, if $M$ is the matrix defined by

$$ (1\ \al\ \cdots\ \al^{d-1})M\ = \ (1\ \g\ \cdots \g^{d-1})
$$
then

$$ (\al_j^k)M=(\g_j^k)
$$
where $\al_j$ ranges over all the conjugates of $\al$ ($j=1,...,d$)
and $k=0,1,...,d-1$. Taking the determinant of both sides,
and using the Vandermonde determinant formula, we see that

$$ \det(M)^2\ = \ {\p_{i\not=j} (\g_i-\g_j)\over \p_{i\not=j} (\al_i-\al_j)}\ = \
{\rm Norm}_{T/S}
\left({f'(\gamma)\over f'(\al)}\right)
$$

Thus the map $\langle f'(\al)\rangle \mapsto  [\det(\al_j^k)]^2$
is a well defined isomorphism from $\langle \cD\rangle$ to $\d$,
and this proves Lemma \ref{different}.

\v

\begin{lemma}\label{DK}
Let $\pi:Y\ra B$ be a finite flat separable morphism of smooth
quasi-projective varieties. Then $\cD_{Y/B}=K_{Y/B}^{-1}$. Thus if
$L\ra Y$ is a line bundle, we have a canonical isomorphism:

\be
\det (\pi_*L)^{ 2} \ \ra \
\langle L\rangle^2\otimes \langle K_{Y/B}^{-1}\rangle\
\ee
Replacing $L$ by $L^k$ yields an isomorphism

\be
\m_k(L,Y,B): \det (\pi_*L^k)^{ 2}\
\ra\ \langle L\rangle^{2k}\otimes \langle K_{Y/B}^{-1}\rangle
\ee

\end{lemma}

{\it Proof.} As before, we may assume $B={\rm spec } (S)$ and $Y={\rm spec}(T)$.
We have an exact sequence \cite{H}

$$ 0\ \ra \ \pi^*\O_{B/\C}\ \ra \ \O_{Y/\C}\ \ra \ \O_{Y/B}\ \ra 0,
$$
where injectivity of the second arrow follows from the fact that
it is a full rank morphism of two vector bundles.
According to (\ref{short}),  $\O_{Y/B}=T/(f'(\al))$.  Thus $H^0(Y,\O_{Y/\C})$ has
a basis $\o_1,...,\o_m$ with the property: $f'(\al)\o_1,\o_2,...,\o_m$ is
a basis of $H^0(Y,\pi^*\O_{B/\C})$. This shows that $K^{-1}_{Y/B}$ is
principal, and generated by $f'(\al)$. Since $f'(\al)$ also
generates $\cD_{Y/B}$, Lemma \ref{DK} is proved.

\section{Arbitrary relative dimension}

\setcounter{equation}{0}
\v
In this section we complete the proof of Theorem \ref{sharp1}.
First we define the isomorphism $\m_k$ of (\ref{phik}).
To do this, we choose sections $s_1,...,s_N$ of $L$ which are in general
position.

The adjunction formula says that
$(K_{X/B}L^n)|_{X_{n-1}}=K_{X_{n-1}/B}(L|_{X_{n-1}})^{n-1}$.
Thus
we have an isomorphism:

\be \ti\iota_{s_N}:\langle L,...L\rangle_{X/B}^{2k}\langle K_{X/B}L^n,...,L\rangle_{X/B}^{-1}
\ \ra \
\langle L,...L\rangle_{X_{n-1}/B}^{2k}\langle K_{X_{n-1}/B}L^{n-1},...,L\rangle_{X_{n-1}/B}^{-1}
\ee

where, on the left side,
there are $n+1$ terms in each pairing and on the right
side there are $n$ terms.
Continuing in this fashion, we obtain an isomorphism:
\be
\ti\iota(s_1,...,s_N)=\ti\iota_{s_1}\circ\cdots\circ\ti\iota_{s_N}:\langle L,...L
\rangle_{X/B}^{2k}\langle
K_{X/B}L^n,...,L\rangle_{X/B}^{-1} \ra \
\langle L_Y^k\rangle^2_{Y/B}\langle K^{-1}_{Y/B}\rangle_{Y/B}
\ee

where $Y=X_0$.

\v

Define the isomorphism
$\ti\kappa(s_1,...,s_n)=\kappa_{s_1}\circ \D\kappa_{s_2}\circ\cdots \circ
\D^{(n-1)}\kappa_{s_n}$. Then

\be \ti\kappa(s_1,...,s_n): \D^{(n)}\det(\pi_*L^k)\ \ra \ \det(\pi_*L^k_{X_0})
\ee

\v

Finally, define

\be
\m_k(L,X,B)(s_1,...,s_n)\ = \ \ti\iota(s_1,...,s_n)^{-1}\circ \mu_k(L|_{X_0},X_0,B)\circ
\ti\kappa^2(s_1,...,s_n)
\ee

\v
As in the proof of Theorem \ref{sharp}, the proof of
Theorem \ref{sharp1} follows from the following :

\begin{lemma}\label{ind1}
The isomorphism $\m_k(L,X,B)(s_1,...,s_n)$ is independent of the choice
of generic sections $s_1,...,s_n$.
\end{lemma}

The first step is to show that $\m_k$ doesn't change
if the sections $s_1,...,s_n$ are permuted. The
proof is similar to the first step in the proof
of Lemma \ref{ind}, so we omit it.

\v

As before, we are reduced to proving Lemma \ref{ind1} in
the case $n=1$:
We shorten the notation by
writing $\m(s)= \m_k(L,X,B)(s)$. Thus

$$\m(s): [\D\det\ \pi_*(L^k)]^2\ \ra \
\langle L,L\rangle^{2k}\langle K_{X/B}L,L\rangle^{-1}\ = \
\langle L^{2k-1}K_{X/B}^{-1},L\rangle
$$

{\bf Step 1.} If $s$ is a generic section, then for every $t\in\C^\times$ we have
\be\label{claim}\m(s)=\m(ts)
\ee

 To
prove this, we first recall the definition of $\m$:

\be\label{(1)}\m(s)\ = \ \tilde\iota(s)^{-1}\circ \m(L|_{X_0},X_0,B)\circ\tilde \kappa(s)^2
\ee

Here $X_0=\{s=0\}$.
We claim that there is an integer $p$ with the property

\be\label{(2)} \tilde \kappa^2(ts)\ = \ t^p\tilde\kappa^2(s)\ \ \ {\rm and}\ \ \
\tilde\iota(ts)\ = \ t^p\tilde\iota(s)\ .
\ee

If we prove (\ref{(2)}) then (\ref{(1)}) implies $\m(s)=\m(ts)$. Thus we need only prove (\ref{(2)}).

\v
To ease the notation, we shall write $\iota$ for $\ti \iota$ and $\kappa$ for $\ti\kappa$.
First we examine $\kappa(ts)$. Recall that $\kappa(s): \D\det(\pi_*L^k)\ra\det(\pi_*L^k_{X_0})$
is defined by

$$ \kappa(s)\left({\det(sF,H)\over \det(F)}\right)\ = \ \det(H)
$$
where $F$ is an ordered basis of $H^0(X,L^{k-1})$ and $H\sub H^0(X,L^k)$ is an ordered
set such that $(sF,H)$ is a basis of $H^0(X,L^k)$.
Replacing $s$ by $ts$ we have

$$ \kappa(s)\left({\det(sF,H)\over \det(F)}\right)\ =
\kappa(ts)\left({\det(tsF,H)\over \det(F)}\right)\ = \ \kappa(ts)t^{|F|}\left({\det(sF,H)\over \det(F)}\right)
$$

Thus $\kappa(s)=t^{|F|}\kappa(ts)$.
\v

Next we calculate $\iota(ts)$. Recall that

$$\iota(s): \langle L^{2k-1}K_{X/B}^{-1},L\rangle \ \ra \
\langle L^{2k}K^{-1}_{X_0/B}\rangle
$$
In other words,

$$\iota(s): \langle L^{1-2k}K_{X/B},L\rangle^{-1} \ \ra \
\langle L^{-2k}K_{X_0/B}\rangle^{-1}
$$
The definition
is given as follows:

\be\label{(3)} \iota(s)(\langle\o,s\rangle^{-1})\ = \ \langle{\rm Ad}(s)(\o)\rangle^{-1}
\ee

where ${\rm Ad}(s):  L^{1-2k}K_{X/B}^{-1}|_{X_0}\ra L^{-2k}K_{X_0/B}$ is
the isomorphism given  by the adjunction formula, and is characterized
by the formula:

$$  {\o\over s}\ = \ {df\over f}\wedge  {\rm Ad}(s)(\o)
$$

where $f$ is any local defining equation of $X_0\sub X$.
\v
Thus, replacing $s$ by $ts$ we see that

$$ {\rm Ad}(ts)(\o)\ = \ t^{-1}{\rm Ad}(s)(\o)
$$
so, if $\deg(L)$ is the degree of $L$ on a fiber of $X\ra B$, we have
$$ \langle{\rm Ad}(ts)(\o)\rangle\ = \ t^{-\deg(L)}\langle{\rm Ad}(s)(\o)\rangle
$$

using properties of the Deligne pairing. Replacing $s$ by $ts$ in (\ref{(3)}) we get

\be\label{(4)} \iota(ts)(\langle\o,ts\rangle^{-1})\ = \ \langle{\rm Ad}(ts)(\o)\rangle^{-1}\ = \
t^{\deg(L)}\langle{\rm Ad}(s)(\o)\rangle^{-1}\ .
\ee
On the other hand, the properties of the Deligne pairing imply

$$ \langle\o,ts\rangle^{-1}\ = \ t^m\langle\o,s\rangle^{-1}\ ,
$$

where $m$ is the degree of $L^{2k-1}K^{-1}$ on a fiber of $X\ra B$.
Thus we conclude

$$ \kappa^2(s)=t^{2|F|}\kappa^2(ts)\ \ \ {\rm and} \ \ \
\iota(s)=t^{m-\deg(L)}\iota(ts)\ .
$$

To prove (\ref{(2)}) we must  show that $2|F|=m-1$, that is, we must show

\be\label{(5)} 2\dim(H^0(X,L^{k-1})\ = \ \deg(L^{2k-1}K^{-1})-\deg(L)\ = \
\deg(L^{2(k-1)}K^{-1})
\ee

The Riemann-Roch formula tells us that $2\dim(H^0(X,L^k))=\deg(L^{2k}K^{-1})$
for $k$ sufficiently large. This now proves (\ref{claim}) and
completes Step 1.

\v

To describe Step 2, we
first recall that

$$\m^{-1}(s): \langle L^{2k-1}K_{X/B}^{-1},L\rangle\ \ra \
[\D\det\ \pi_*(L^k)]^2
$$

To prove Lemma \ref{ind1} we must show that if $s,s'$ are generic sections,
then $\m(s)=\m(s')$. To do this, we connect $s$ to $s'$ by a line $s_t=(u+t)s'=s+ts'$,
and study $\m(s_t)$ as $t$ varies.
But $\m(s_t)$ is only defined if $Y_t=\{s_t=0\}$ is smooth
and flat over $B$.
Bertini's theorem  says that $Y_t$ is smooth
and flat over $B$ for {\it all but finitely many $t$}. On
the other hand, by Proposition \ref{Bertini 2}, if $s''$ is a generic section, then
$\{s+ts''=0\}$ and $\{s'+ts''=0\}$ are smooth and flat for all $t\in\C$
(possibly after shrinking the base $B$).
Thus we may assume that $Y_t$ is smooth and flat over $B$ for all $t\in\C$
so that $\m(s_t)$ is well defined for {\it all} $t\in\C$.

\v\v
Thus we let  $s,s'\in H^0(X,\pi_*L)$ be generic sections
with the property: $Y_t=\{s+ts'=0\}$ is smooth and flat
over $B$ for all $t\in\C$.
Let
$\eta\in H^0(X,\pi_*(L^{2k-1}K^{-1})) $
be such that
$(\eta)\cap (s)=(\eta)\cap (s')=\emptyset$.
To prove Lemma \ref{ind1}
we must show that

\be\label{to show}
 \m^{-1}(s)(\langle \eta, s\rangle)\ = \ \m^{-1}(s')(\langle \eta, s\rangle)
\ee

\v

We write $s=us'$ where $u$ is a rational
function such that $(u)\cap (\eta)=\emptyset$.
Thus

\be\label{s and sprime}
 \langle \eta,s\rangle\ = \ u\big[(\eta)\big]\langle \eta,s'\rangle
\ee
where, for $b\in B$,
\be u\big[(\eta)\big](b)\ = \ \p_{q\in (\eta)\cap X_b} u(q)
\ee
\v

{\bf Step 2.}
Let
$m$ be the degree of $L^{2k-1}K^{-1}$ on a
fiber of $X\ra B$.
Then we claim that

\be\label{limit1}
\label{note} \m^{-1}(s')(\langle\eta,s'\rangle)=\lim_{t\to\i}t^{-m}\m^{-1}(s_t)(\langle\eta,s_t\rangle)\ .
\ee
To prove  (\ref{limit1}), first recall that

$$\m^{-1}(s): \langle L^{2k-1}K_{X/B}^{-1},L\rangle\ \ra \
[\D\det\ \pi_*(L^k)]^2
$$

We must show that if $s,s'$ are generic, then
$\m(s)=\m(s')$, which is equivalent to showing
that $\m^{-1}(s)=\m^{-1}(s')$.

\v

We now return to the proof of (\ref{limit1}). Recall that $s_t=ts'+s$. Relation
(\ref{claim}) implies $\m(s_t)=\m(s'+{1\over t}s)$.
We can thus rewrite (\ref{limit1})
as follows:

$$ \lim_{t\to 0}\m^{-1}(s'+ts)\langle(\eta, s'+{ t}s\rangle)\ = \
\m^{-1}(s')(\langle\eta,s'\rangle)
$$

Thus, to prove (\ref{limit1}), we must show that for fixed $\eta$, the map

\be\label{cont}
F(t) \  =\ \m^{-1}(s'+ts)(\langle\eta,s'+ts\rangle)
\ee

is a continuous function of $t$ and a neighborhood of $t=0$.
Here $s,s'\in H^0(B,\pi_*L)$ and
$\eta\in H^0(B,\pi_*(L^{2k-1}K_{X/B}^{-1}))$. To do this, we let $\ti X=X\times\C$,
$\ti B=X\times \C$ and we let $\ti\pi:\ti X\ra \ti B$ be
the map $\ti \pi(x,t)=(\pi(x),t)$. We
let $p:\ti X\ra X$ be the projection map
and $\ti L=p^*L$ so that $\ti L\ra \ti X$ is
a line bundle.
\v

If $\ti s\in H^0(\ti B,\ti\pi_*\ti L)$
and
$\ti \eta\in
H^0(\ti B,
\ti \pi_*(\ti L^{2k-1}K_{\ti X/\ti B}^{-1}))$ then
$\ti \m^{-1}(\ti s)(\langle\ti \eta,\ti s\rangle)$ is a section of
$[\D\det\ti\pi_*(\ti L^k)]^2
$, that is,
$$\ti \m^{-1}(\ti s)(\langle\ti \eta,\ti s\rangle):B\times\C \ra [\D\det\ti\pi_*(\ti L^k)]^2
$$

Now fix $t_0\in\C$ and let
$\sigma_{t_0}: B\ra B\times \C$ be the
map $\sigma_{t_0}(b)=(b,t_0)$. Then
$$\sigma_{t_0}^*[\D\det\ti\pi_*(\ti L^k)]^2=
[\D\det\pi_*( L^k)]^2$$

With these preliminaries, we return to
the proof of the continuity of
$F$. Fix $s,s'$
as before and let $\ti s: X\times\C\ra L$
be the map $\ti s(x)=s'(x)+ts(x)$.
Thus $\ti s\in  H^0(\ti B,\ti\pi_*\ti L)$,
in other words, $\ti s$ is a section of
$\ti L\ra \ti X$.

\v

It follows from the definition of $\mu$ and $\ti\mu$ that for each $t_0\in\C$
and $b\in B$ we have

$$\m^{-1}(s'+t_0s)(\langle\eta,s'+ts\rangle)(b)\ = \ \sigma_{t_0}^*[\ti \m^{-1}(\ti s)
(\langle\ti \eta,\ti s\rangle)(b,t_0)]
$$

Now the right side of this
equation is manifestly continuous
(in fact, analytic) in the
$t_0$ variable. This shows that $F(t)$ is an analytic
function of $t$, and thus continuous at $t=0$.
This completes Step 2.

\v\v
Now we can finish the proof of Lemma \ref{ind1}:
for $t\in\C$ with $(\eta)\cap (s_t)=\emptyset$, we
define $\rho_\eta(b;t)\in\C^\times$ by the formula

\be\label{rho0}
\m(s)\m(s_t)^{-1}(\langle\eta,s_t\rangle)\ = \ \r_\eta(b;t)\langle \eta, s\rangle
\ee
\v
We define $\r_\eta(b;t)$ for an arbitrary number $t\in\C$ as follows:
Suppose $t_0\in\C$ is such that $(\eta)\cap (s_{t_0})\not=\emptyset$.
Then choose $\eta^*$ a generic section of $H^0(B,\pi_*(L^{2k-1}K^{-1}))$.
Then

$$ \r_{\eta^*}(b;t){\eta^*\over \eta}((s))\cdot{\eta\over \eta^*}((s_t))\ = \
\r_{\eta}(b;t)
$$

The left side is well defined in a neighborhood of $t=t_0$
and this shows that $\r_\eta(b;t)$ extends to an analytic
function on all of $t\in\C$ which vanishes precisely when
$(\eta)\cap (s_t)\not=\emptyset$, that is, precisely
when $u(q)=-t$ for some $q\in (\eta)$. Moreover, according to
(\ref{note}),  $t^{-m}\r_\eta(b;t)$ has a finite limit as
$t\ra \i$. This shows that $\r_\eta(b;t)$ is a polynomial
of degree $m$. Thus we have

$$  \r_\eta(b;t)\ = \ \al(b)\p_{q\in(\eta)} (t+u(q))
$$

with $\al(b)\not=0$.
Since
$\rho_\eta(b;0)=1$, we must have $\al(b)=u\big[(\eta)\big]^{-1}$.
Now (\ref{rho0}) implies
\be\label{rho01}
t^{-m}\m(s_t)^{-1}(\langle\eta,s_t\rangle)\ = \ t^{-m}\r_\eta(b;t)\m^{-1}(s)\langle \eta, s\rangle
\ee
Taking the limit as $t\to\i$, and applying (\ref{limit1}), we obtain

\be \m^{-1}(s')(\langle \eta, s'\rangle)\ = \ u[(\eta)]^{-1}\m^{-1}(s)\langle\eta,s\rangle
\ee

Thus (\ref{to show}) follows from (\ref{s and sprime}).

\section{Asymptotics of the Mabuchi K-energy}

\v
In this section we prove Corollary 2 (the proof of
Corollary 3 is exactly the same, so we omit it).
First we recall some properties of the Deligne metric: Let $\pi:X\ra B$
be a flat projective morphism between smooth quasi-projective
varieties. Let $n$ be the relative dimension of $\pi$, and for $0\leq j\leq n$,
let $h_j$ be a smooth metric on $L_j$. Let $\langle h_0,...,h_n\rangle$ be
the Deligne metric on $\langle L_0,...,L_n\rangle$ as defined by
Deligne \cite{Del} and
Zhang~\cite{Z96}. This metric
is continuous by the result of Moriwaki \cite{M} and it satisfies
the following change of metric property \cite{PS03}:

\be \label{small change}
\langle h_0,...,h_{n-1},h_ne^{-\phi}\rangle \ = \
\langle h_0,...,h_n\rangle e^{-\Phi}
\ee
where
\be\label{2.6}  \Phi \ = \ \I_{X/B} \    \phi\cdot
\p_{k <n} c'_1(h_k)
\ee

and $c_1'(h_k)=-{i\over 2\pi}\ddb\log h_k$.

\v

Next we recall the formula for $\eta$:

\be\label{PS}
\eta(L,X)\ = \ \langle L,...,L\rangle^\m
\otimes
\langle K,L,...,L\rangle^{(n+1)}
\ee

If $h$ is a metric on $L$ with positive curvature $\o$,
then ${1\over \o^n}$ is a smooth metric on
$K$ and thus we obtain a Deligne metric

\be
\eta(h)=\langle h,...,h\rangle^\m\otimes \langle {1\over \o^n},h,...,h\rangle^{(n+1)}
\ee

which is smooth on $X'\sub X$, the union of all the smooth fibers.
\v
The key property of the metric $\eta(h)$ is the following transformation
rule \cite{PS02, PS03}: Let $\o'=\o+{\sqrt{-1}\over 2\pi}\phi$ be a \K metric in the same
\K class as $\o$.  Then

\be\label{K change}
 \eta(he^{-\phi})\ = \ \eta(h)e^{-\ti\n(\o,\o')}
\ee
where $\ti\n(\o,\o')=d(n+1)\n(\o,\o')$,
$d$ is the
degree of $L$ on a smooth fiber,
 and $\n(\o,\o')$ is
the K-energy of $\o'$ with respect to $\o$. In fact, one
may use (\ref{K change}) to define $\n(\o,\o')$. Since
the map $\r(t): (L_1,h_1e^{-\phi})\ra (L_t,h_t)$ is
an isometry (this is the definition of $\phi$) we conclude
that for any non-zero section $s_1$ of $\eta(L,X)$,

\be
\n(t)\ = \ \n(\o_1,\r(t)^*\o_t)\ = \ {1\over d(n+1)}\log{\|\r(t)s_1\|^2\over \|s_1\|^2}
\ee
where $\|\cdot\|$ is the metric defined by $\eta(h)$.
\v

{\it Proof of Corollary 2.}
 Let

\be
\eta^*(h)=\langle h,...,h\rangle^\m\otimes \langle {\pi^*(dt\wedge d\bar t)\over \o^{n+1}},h,...,h\rangle^{(n+1)}
\ee

Then $\eta^*(h)$ is a continuous metric on all of $\eta(L,X)$.
Moreover, (\ref{small change}) implies

\be \eta^*(h)\ = \ \eta(h)e^{-d(n+1)\psi}
\ee

Combining this with (\ref{K change}) we see that if $s_1$ is
a non-zero element of  $L_1$, then

\be {1\over d(n+1)}\cdot\log{\|\r(t)s_1\|^2_*\over \|s_1\|^2_*}\ = \
\n(t)\ - \ \psi(t) + \psi(1)
\ee
where $\|\cdot\|_*$ is the metric defined by $\eta^*(h)$.
\v
To complete the proof of Corollary 2, we need the following:

\begin{lemma}
\be  \log{\|\r(t)s_1\|_*\over \|s_1\|_*}\ =\ F(T)\log |t|\ + \ \e(t)
\ee
where $\e(t)$ is a continuous function in a neighborhood of $t=0$.
\end{lemma}
{\it Proof.} Let $s:\C\ra\eta(L,X)$ be a nowhere vanishing section. Define, for $t\not=0$
and $z\in\C$, the function $f(t,z)$ by the formula

\be\label{rho} \r(t)(s(z))\ = \ f(t,z)s(tz)
\ee

Then Theorem \ref{main} together with \cite{PT06}
implies $f(t,0)=t^{F(T)}$. Applying $\r(t')$ to both sides of
(\ref{rho}) we get

\be\label{rho1}
\r(t')\r(t)(s(z))\ = \ f(t',tz)f(t,z)s(t'tz)\ = \ f(t't,z)s(t'tz)
\ee
Thus

\be
\label{g}
f(t,z)\ = \ {f(tz,1)\over f(z,1)}\ = \ {g(tz)\over g(z)}
\ee

for all $z\not=0$, where $g(t)=f(t,1)$. Since the right side of (\ref{g})
approaches $t^{F(T)}$ as $z$ approaches zero, we see that $g$
does not have an essential singularity at the origin.
Define an integer
\be
q\in\Z
\ee
such that $h(t)=t^{-q}g(t)$ is
non-zero and holomorphic in a neighborhood of $t=0$. Then
$f(t,z)=t^q h(tz)/h(z)$. This shows that
\be
q=F(T)
\ee
and
\be
\log{\|\r(t)s_1\|_*\over \|s_1\|_*}\ =\ F(T)\log|t|\ +
\e(t)
\ee
where
\be
\e(t)\ = \ \log|h(t)/h(1)|\ + \ \log \|s(t)\|_*\ - \log\|s_1\|_*
\ee

By Moriwaki \cite{M}, the term $\log \|s(t)\|_*$ is continuous.
Since $h(t)$ is holomorphic and non-vanishing near $t=0$, we
conclude that $\e(t)$ is continuous.

\v
{\it Remark 6.} It has been pointed out to us by Shou-Wu Zhang
that $q$ can be viewed as a non-archimedian Mabuchi functional
on the space of test configurations.

\enddocument